# ESTIMATION OF THE COVARIANCE MATRIX OF RANDOM EFFECTS IN LONGITUDINAL STUDIES

By Yan Sun, Wenyang Zhang and Howell Tong

*Shanghai University of Finance and Economics, University of Kent and London School of Economics*

Longitudinal studies are often conducted to explore the *cohort* and *age* effects in many scientific areas. The within cluster correlation structure plays a very important role in longitudinal data analysis. This is because not only can an estimator be improved by incorporating the within cluster correlation structure into the estimation procedure, but also the within cluster correlation structure can sometimes provide valuable insights in practical problems. For example, it can reveal the correlation strengths among the impacts of various factors. Motivated by data typified by a set from Bangladesh pertinent to the use of contraceptives, we propose a random effect varying-coefficient model, and an estimation procedure for the within cluster correlation structure of the proposed model. The estimation procedure is optimization-free and the proposed estimators enjoy asymptotic normality under mild conditions. Simulations suggest that the proposed estimation is practicable for finite samples and resistent against mild forms of model misspecification. Finally, we analyze the data mentioned above with the new random effect varying-coefficient model together with the proposed estimation procedure, which reveals some interesting sociological dynamics.

**1. Introduction.**

1.1. *Theoretical background.* Longitudinal data analysis has attracted considerable attention from statisticians recently. The methodology for parametric-based longitudinal data analysis is quite mature; see, for example, Diggle, Heagerty, Liang and Zeger [5] and the references therein. The situation with nonparametric-based longitudinal data analysis is quite different.

One of the main difficulties is how to incorporate the within cluster correlation structure into the estimation procedure. For nonparametric longitudinal regression, see Zeger and Diggle [30], Hoover, Rice, Wu and Yang [15], Fan and Zhang [10], Lin and Carroll [19], Wu and Zhang [28], Fan and Li [8], Qu and Li [24] and others. He, Fung and Zhu [13] investigate the robust estimation in generalized partial linear models for longitudinal data. Lin and Carroll [19] recommend that we ignore the within cluster correlation when kernel smoothing is employed. Welsh, Lin and Carroll [26] investigate the possibility of using weighted least squares based on the within cluster correlation structure when spline smoothing is used. They suggest that the weighted least squares estimator works better than the estimator based on working independence when spline smoothing is used. Wang [25] provides an innovative kernel smoothing and demonstrates that, when the true correlation is available, her estimator is more efficient than the most efficient estimator that is obtained by adopting a working independence approach.

In longitudinal data analysis, whether parametric or nonparametric, within cluster correlation structure can be used to improve the efficiency of the estimation. However, the within cluster correlation structure is typically unknown in reality. In this paper, we investigate systematically the estimation of the within cluster correlation structure.

1.2. *Practical meaning.* The within cluster correlation structure can lead to not only important improvement of the estimation but also some practical insights. As we shall see in the analysis of the Bangladesh data, the estimated within cluster correlation structure actually sheds interesting light on the impacts among factors.

The Bangladesh data set is from the Bangladesh Demographic and Health Survey 1996–1997. This survey follows a two-stage sample design in which clusters were selected at the first stage, and women were sampled from these clusters at the second stage. The clusters correspond to villages in rural areas and neighborhoods in urban areas, and may loosely be termed communities. What is of interest is how the factors which are commonly found to be associated with fertility in Bangladesh affect the first birth interval, and how strongly correlated are the effects of these factors. The selected factors are (1) age at first marriage; (2) woman's level of education; (3) type of region of residence; (4) woman's religion; (5) year of marriage; (6) administrative area. Among these factors, type of region of residence and administrative area pertain to cluster levels and are called cluster-level variables, and the rest are called individual-level variables.

We use $y$ to denote the length of the first birth interval, $Z$ the vector of individual-level variables, $\xi$ the vector of cluster-level variables and **e** the random effect. For $j = 1, \ldots, n_i$, $i = 1, \ldots, m$, let $y_{ij}$ and $Z_{ij}$ be the $j$th observation of $y$ and $Z$ in the $i$th cluster, $\xi_i$ the observation of $\xi$ at the $i$th cluster, and $\mathbf{e}_i$ the random effect of the $i$th cluster, which is unobservable.

When examining the effects of year of marriage and other factors on the length of the first birth interval, it is necessary to take into account clustering of responses for women in the same community. This is because the first birth intervals for women in the same cluster may be correlated due to unobserved cluster-level factors such as cultural norms and access to family planning programes. The usual way to incorporate unobservable variables in a statistical model is via random effects. This leads to the multilevel model

$$(1.1) \qquad y_{ij} = Z_{ij}^{\mathrm{T}}(\mathbf{a}_1 + \mathbf{e}_i) + \xi_i^{\mathrm{T}}\mathbf{a}_2 + \varepsilon_{ij}.$$

The coefficient $\mathbf{a}_1 + \mathbf{e}$ can be regarded as the impact of $Z$ on $y$, which is random across the clusters. The correlation matrix of $\mathbf{e}$ can reveal how strongly correlated are the impacts of the components of $Z$ on $y$. In this paper, we propose an estimation procedure for the covariance matrix of $\mathbf{e}$.

Let

$$X_{ij} = (Z_{ij}^{\mathrm{T}}, \xi_i^{\mathrm{T}})^{\mathrm{T}}, \qquad \mathbf{a} = (\mathbf{a}_1^{\mathrm{T}}, \mathbf{a}_2^{\mathrm{T}})^{\mathrm{T}}.$$

Equation (1.1) can be written as

$$(1.2) \qquad y_{ij} = X_{ij}^{\mathrm{T}}\mathbf{a} + Z_{ij}^{\mathrm{T}}\mathbf{e}_i + \varepsilon_{ij}.$$

Model (1.2) assumes that the impacts of the factors on the length of the first birth interval are time-invariant, which may not be plausible because as a society Bangladesh is changing with time. So, it is more realistic to assume that the impacts can vary with time. This leads to the following random effect varying-coefficient model:

$$(1.3) \qquad y_{ij} = X_{ij}^{\mathrm{T}}\mathbf{a}(U_{ij}) + Z_{ij}^{\mathrm{T}}\mathbf{e}_i + \varepsilon_{ij}, \qquad j = 1, \ldots, n_i, i = 1, \ldots, m,$$

where $\varepsilon_{ij}$, $j = 1, \ldots, n_i$, $i = 1, \ldots, m$, are measurement errors, which we assume to be i.i.d. with $E(\varepsilon_{ij}) = 0$ and $\mathrm{var}(\varepsilon_{ij}) = \sigma^2$. Here $\mathbf{e}_i$, $i = 1, \ldots, m$, are random effects across the clusters, which we assume to be i.i.d. with $E(\mathbf{e}_i) = 0$ and $\mathrm{cov}(\mathbf{e}_i) = \Sigma$. Further, $\mathbf{e}_i$ is independent of $\varepsilon_{ij}$; $X_{ij}$ is a $p$-dimensional covariate and $Z_{ij}$ is a $q$-dimensional covariate associated with random effects. We assume that $n_i < N < \infty$. $(U_{ij}, X_{ij}, Z_{ij})$ are i.i.d. and independent of $\mathbf{e}_i$ and $\varepsilon_{ij}$.

The within cluster correlation structure in (1.3) has been used extensively in the literature to model the cluster effect. See, for example, Laird and Ware [18], Jennrich and Schluchter [17], Longford [21], Zeger, Liang and Albert [31], Lindstrom and Bates [20], Hedeker and Gibbons [14] and others. The focus in the above cited was on the estimation relevant to the regressive coefficients in parametric models. Wu and Liang [27] proposed an interesting backfitting estimation procedure to estimate the functional coefficients in a time-varying-coefficient mixed effects model.

Varying-coefficient models are useful when exploring dynamic systems. There is a growing literature addressing this kind of model, which includes Hastie and Tibshirani [12], Xia and Li [29], Cai, Fan and Li [1], Zhang, Lee and Song [32], Fan, Yao and Cai [9] and the references therein.

There is some literature discussing how to estimate the within cluster covariance matrix for the parametric setting. The commonly used method is restricted maximum likelihood estimation (REML); see Laird and Ware [18]. While REML is theoretically appealing, the optimization involved can be very difficult. Recently, Fan, Huang and Li [7] proposed an innovative semiparametric estimation procedure for the covariance structure in longitudinal studies. They investigated both the quasi-likelihood approach and the minimum generalized variance approach. In this paper, focusing on model (1.3), we take a different approach to studying the estimation of $\Sigma$ and $\sigma^2$ that includes both the methodology and relevant asymptotics. The estimators proposed in this paper have closed forms, so they are easy to implement without any need for optimization. Also, the proposed estimation method does not depend on the likelihood function, so it is robust against likelihood specification. We have established asymptotic normality of the proposed estimators. We have also conducted extensive simulation studies, which show that when the dimension of $\mathbf{e}_i$ is not larger than 3, REML encounters no serious optimization problem, and in such cases REML and our estimators have comparable performance, with perhaps the latter being marginally better. On the other hand, when the dimension of $\mathbf{e}_i$ is larger than 3, we have not been able to find a convergent algorithm for REML. In this paper, we will also show that, when spline smoothing is used, the weighted least squares estimators of the functional coefficients perform much better if we incorporate the within cluster correlation structure estimated by our proposed method instead of assuming working independence.

Our paper is organized as follows. We begin in Section 2 with a description of an estimation procedure for $\Sigma$ and $\sigma^2$. We present the asymptotic properties of the proposed estimators in Section 3 and assess the performance of the model by simulation in Section 4. In Section 5, using the new model and the proposed estimation procedure, we explore how the impacts of the factors on the length of first birth interval in Bangladesh change with time and how strongly correlated are the impacts of the factors on the length of first birth interval.

**2. Estimation procedure.** The procedure first estimates $\mathbf{a}(\cdot)$ based on working independence, then uses the residual to estimate $\sigma^2$, and finally estimates $\Sigma$.

For any given $u$, we use $\mathbf{a}$ and $\dot{\mathbf{a}}$ to denote $\mathbf{a}(u)$ and $d\mathbf{a}(u)/du$, respectively. By Taylor's expansion, we have

$$\mathbf{a}(U_{ij}) \approx \mathbf{a} + \dot{\mathbf{a}}(U_{ij} - u)$$

when $U_{ij}$ is in a small neighborhood of $u$. This leads to the local least squares procedure

$$\text{(2.1)} \qquad \sum_{i=1}^{m}\sum_{j=1}^{n_i}[y_{ij} - X_{ij}^{\text{T}}\{\mathbf{a} + \mathbf{b}(U_{ij} - u)\}]^2 K_h(U_{ij} - u),$$

where $K_h(\cdot) = K(\cdot/h)/h$, $K(\cdot)$ is the kernel function and $h$ is the bandwidth. We minimize (2.1) with respect to $(\mathbf{a}, \mathbf{b})$ to get the minimizer $(\hat{\mathbf{a}}, \hat{\mathbf{b}})$. The estimator of $\mathbf{a}$ is taken to be $\hat{\mathbf{a}}$. By simple calculations, we have

$$\hat{\mathbf{a}} = (I_p, \mathbf{0}_p)\left(\sum_{i=1}^{m}\Lambda_i^{\text{T}}W_i\Lambda_i\right)^{-1}\sum_{i=1}^{m}\Lambda_i^{\text{T}}W_iY_i,$$

where $I_p$ is $p \times p$ identity matrix and $\mathbf{0}_p$ is the $p \times p$ matrix with all entries being 0,

$$\Lambda_i = (\mathbf{X}_i, D_i\mathbf{X}_i), \qquad \mathbf{X}_i = (X_{i1}, \ldots, X_{in_i})^{\text{T}},$$
$$D_i = \text{diag}(U_{i1} - u, \ldots, U_{in_i} - u),$$
$$W_i = \text{diag}(K_h(U_{i1} - u), \ldots, K_h(U_{in_i} - u)), \qquad Y_i = (y_{i1}, \ldots, y_{in_i})^{\text{T}}.$$

Next, we estimate $\sigma^2$. Let $\hat{\mathbf{a}}(U_{ij})$ be $\hat{\mathbf{a}}$ with $u$ being replaced by $U_{ij}$. Let

$$\mathbf{r}_i = (r_{i1}, \ldots, r_{in_i})^{\text{T}}, \qquad r_{ij} = y_{ij} - X_{ij}^{\text{T}}\mathbf{a}(U_{ij}), \qquad \hat{\mathbf{r}}_i = (\hat{r}_{i1}, \ldots, \hat{r}_{in_i})^{\text{T}},$$
$$\hat{r}_{ij} = y_{ij} - X_{ij}^{\text{T}}\hat{\mathbf{a}}(U_{ij}), \qquad \mathbf{Z}_i = (Z_{i1}, \ldots, Z_{in_i})^{\text{T}}, \qquad P_i = \mathbf{Z}_i(\mathbf{Z}_i^{\text{T}}\mathbf{Z}_i)^{-1}\mathbf{Z}_i^{\text{T}}.$$

For each given $i$, based on the residual $\mathbf{r}_i$, we have the following synthetic linear model:

$$\text{(2.2)} \qquad \mathbf{r}_i = \mathbf{Z}_i\mathbf{e}_i + \boldsymbol{\varepsilon}_i, \qquad \boldsymbol{\varepsilon}_i = (\varepsilon_{i1}, \ldots, \varepsilon_{in_i})^{\text{T}}.$$

The residual sum of squares of this linear model,

$$\text{rss}_i = \mathbf{r}_i^{\text{T}}(I_{n_i} - P_i)\mathbf{r}_i,$$

would be the raw material for estimating $\sigma^2$. The synthetic degrees of freedom of $\text{rss}_i$ is $n_i - q$. Let $\text{RSS}_i$ be $\text{rss}_i$ with $\mathbf{r}_i$ replaced by $\hat{\mathbf{r}}_i$. $\text{RSS}_i$ is a natural estimator for $\text{rss}_i$. Pooling all $\text{RSS}_i$, $i = 1, \ldots, m$, together naturally leads to the estimator of $\sigma^2$ as

$$\hat{\sigma}^2 = (n - qm)^{-1}\sum_{i=1}^{m}\text{RSS}_i, \qquad n = \sum_{i=1}^{m}n_i.$$

Finally, we estimate $\Sigma$. From (2.2), we have the least squares estimator of $\mathbf{e}_i$ as

$$\tilde{\mathbf{e}}_i = (\mathbf{Z}_i^{\text{T}}\mathbf{Z}_i)^{-1}\mathbf{Z}_i^{\text{T}}\mathbf{r}_i = \mathbf{e}_i + (\mathbf{Z}_i^{\text{T}}\mathbf{Z}_i)^{-1}\mathbf{Z}_i^{\text{T}}\boldsymbol{\varepsilon}_i,$$

which leads to

$$\sum_{i=1}^{m}\tilde{\mathbf{e}}_i\tilde{\mathbf{e}}_i^{\mathrm{T}} = \sum_{i=1}^{m}\mathbf{e}_i\mathbf{e}_i^{\mathrm{T}} + \sum_{i=1}^{m}(\mathbf{Z}_i^{\mathrm{T}}\mathbf{Z}_i)^{-1}\mathbf{Z}_i^{\mathrm{T}}\boldsymbol{\varepsilon}_i\boldsymbol{\varepsilon}_i^{\mathrm{T}}\mathbf{Z}_i(\mathbf{Z}_i^{\mathrm{T}}\mathbf{Z}_i)^{-1} + \sum_{i=1}^{m}(\mathbf{Z}_i^{\mathrm{T}}\mathbf{Z}_i)^{-1}\mathbf{Z}_i^{\mathrm{T}}\boldsymbol{\varepsilon}_i\mathbf{e}_i^{\mathrm{T}}$$
$$+ \sum_{i=1}^{m}\mathbf{e}_i\boldsymbol{\varepsilon}_i^{\mathrm{T}}\mathbf{Z}_i(\mathbf{Z}_i^{\mathrm{T}}\mathbf{Z}_i)^{-1}.$$

The last two terms are of order $O_P(m^{1/2})$, so they are negligible. This leads to

$$m^{-1}\sum_{i=1}^{m}\mathbf{e}_i\mathbf{e}_i^{\mathrm{T}} \approx m^{-1}\left\{\sum_{i=1}^{m}\tilde{\mathbf{e}}_i\tilde{\mathbf{e}}_i^{\mathrm{T}} - \sum_{i=1}^{m}(\mathbf{Z}_i^{\mathrm{T}}\mathbf{Z}_i)^{-1}\mathbf{Z}_i^{\mathrm{T}}\boldsymbol{\varepsilon}_i\boldsymbol{\varepsilon}_i^{\mathrm{T}}\mathbf{Z}_i(\mathbf{Z}_i^{\mathrm{T}}\mathbf{Z}_i)^{-1}\right\}$$
$$\approx m^{-1}\left\{\sum_{i=1}^{m}\tilde{\mathbf{e}}_i\tilde{\mathbf{e}}_i^{\mathrm{T}} - \sigma^2\sum_{i=1}^{m}(\mathbf{Z}_i^{\mathrm{T}}\mathbf{Z}_i)^{-1}\right\}.$$

So, we use

$$(2.3) \qquad \hat{\Sigma} = m^{-1}\sum_{i=1}^{m}\hat{\mathbf{e}}_i\hat{\mathbf{e}}_i^{\mathrm{T}} - m^{-1}\hat{\sigma}^2\sum_{i=1}^{m}(\mathbf{Z}_i^{\mathrm{T}}\mathbf{Z}_i)^{-1}$$

to estimate $\Sigma$. In (2.3), $\hat{\mathbf{e}}_i$ is $\tilde{\mathbf{e}}_i$ with $\mathbf{r}_i$ replaced by $\hat{\mathbf{r}}_i$.

**3. Asymptotic properties.** For any $q \times q$ symmetric matrix $A$, we use $\mathrm{vech}(A)$ to denote the vector consisting of all elements on and below the diagonal of the matrix $A$, $\mathrm{vec}(M)$ to denote the vector by simply stacking the column vectors of matrix $M$ below one another, and let $c_1 = \lim_{m\to\infty} n/(n-qm)$ and $c_2 = \lim_{m\to\infty} n/m$. Obviously there exists a unique $q^2 \times q(q+1)/2$ matrix $R_q$ such that $\mathrm{vec}(A) = R_q \mathrm{vech}(A)$.

To make the presentation more clear, we introduce the following notation. Set

$$\mu_i = \int t^i K(t)\,dt, \qquad i=0,1,2,3, \qquad \boldsymbol{\eta}_i = (X_{i1}^{\mathrm{T}}\mathbf{a}''(U_{i1}),\ldots,X_{in_i}^{\mathrm{T}}\mathbf{a}''(U_{in_i}))^{\mathrm{T}},$$

$$b = (n-qm)^{-1}\sum_{i=1}^{m}\boldsymbol{\eta}_i^{\mathrm{T}}(I_{n_i}-P_i)\boldsymbol{\eta}_i, \qquad B_1 = m^{-1}\sum_{i=1}^{m}(\mathbf{Z}_i^{\mathrm{T}}\mathbf{Z}_i)^{-1},$$

$$B_2 = m^{-1}\sum_{i=1}^{m}(\mathbf{Z}_i^{\mathrm{T}}\mathbf{Z}_i)^{-1}\mathbf{Z}_i^{\mathrm{T}}\boldsymbol{\eta}_i\boldsymbol{\eta}_i^{\mathrm{T}}\mathbf{Z}_i(\mathbf{Z}_i^{\mathrm{T}}\mathbf{Z}_i)^{-1}.$$

Further, we write

$$\Gamma = \lim_{m\to\infty} m^{-1}\sum_{i=1}^{m}E[(\mathbf{Z}_i^{\mathrm{T}}\mathbf{Z}_i)^{-1}],$$

$$\Delta_2 = \lim_{m \to \infty} m^{-1} \sum_{i=1}^{m} E[\text{vec}\{(\mathbf{Z}_i^{\text{T}} \mathbf{Z}_i)^{-1}\} \text{vec}^{\text{T}}\{(\mathbf{Z}_i^{\text{T}} \mathbf{Z}_i)^{-1}\}],$$

$$\Delta_3 = \lim_{m \to \infty} m^{-1} \sum_{i=1}^{m} \sum_{j=1}^{n_i} E[\text{vec}\{(\mathbf{Z}_i^{\text{T}} \mathbf{Z}_i)^{-1} Z_{ij} Z_{ij}^{\text{T}} (\mathbf{Z}_i^{\text{T}} \mathbf{Z}_i)^{-1}\}.$$

$$\times \text{vec}^{\text{T}}\{(\mathbf{Z}_i^{\text{T}} \mathbf{Z}_i)^{-1} Z_{ij} Z_{ij}^{\text{T}} (\mathbf{Z}_i^{\text{T}} \mathbf{Z}_i)^{-1}\}].$$

$$\gamma = \lim_{m \to \infty} (n - qm)^{-1} \sum_{i=1}^{m} \sum_{j=1}^{n_i} E[Z_{ij}^{\text{T}} (\mathbf{Z}_i^{\text{T}} \mathbf{Z}_i)^{-1} Z_{ij}]^2 - c_1 q/c_2 + 1.$$

As the data are unbalanced, that is, different subjects have different numbers of observations, the above expectations take no simple forms. Moreover, let

$$\Delta_1 = \begin{pmatrix} \Sigma \otimes \Gamma_{(1)} + \Gamma \otimes \Sigma_{(1)} \\ \vdots \\ \Sigma \otimes \Gamma_{(q)} + \Gamma \otimes \Sigma_{(q)} \end{pmatrix},$$

where $\Gamma_{(r)}, \Sigma_{(r)}$ $(r = 1, \ldots, q)$ denote the $r$th row of $\Sigma$, $\Gamma$, respectively, and $\otimes$ is the Kronecker product.

THEOREM 1. *Under the technical conditions in the Appendix, we have*

$$n^{1/2} \left\{ \text{vech}(\hat{\Sigma} - \Sigma) - \frac{1}{4} h^4 \left( \frac{\mu_1 \mu_3 - \mu_2^2}{\mu_0 \mu_2 - \mu_1^2} \right)^2 \{\text{vech}(B_2) - \text{vech}(B_1)b\} \right\}$$

$$\xrightarrow{D} N(\mathbf{0}, (R_q^{\text{T}} R_q)^{-1} R_q^{\text{T}} \Delta R_q (R_q^{\text{T}} R_q)^{-1} c_2),$$

*where*

$$\Delta = E\{\mathbf{e}_1 \mathbf{e}_1^{\text{T}} \otimes \mathbf{e}_1 \mathbf{e}_1^{\text{T}}\} - \text{vec}(\Sigma) \text{vec}^{\text{T}}(\Sigma) + \sigma^2 \{\Sigma \otimes \Gamma + \Gamma \otimes \Sigma + \Delta_1\}$$
$$+ 2\sigma^4 \{\Delta_2 - \Delta_3 + c_1/c_2(1 + \gamma)\Gamma\} + \text{var}(\varepsilon_{11}^2)\{\Delta_3 + c_1/c_2 \gamma \Gamma\}.$$

It is clear from Theorem 1 that the estimator $\hat{\Sigma}$ would achieve root-$n$ convergence rate if the bandwidth is properly selected, say the bandwidth $h$ is taken to be $O(n^{-1/8})$. For the estimation of the regression function based on the within cluster correlation structure, Welsh, Lin and Carroll [26] suggest that the spline-based weighted least squares estimation with the right weight would have smaller variance than the working independence approach, but they do not take the bias into consideration. Bias and variance are equally important when assessing the goodness of an estimator. To appreciate both bias and variance, it is better to use the mean squared error as a criterion to assess the accuracy of an estimator. It is not clear whether Welsh, Lin and Carroll's estimator is more efficient than the working independence one in terms of the mean squared error. How to construct a good estimator of the regression function is very important and interesting, but it lies beyond the scope of this paper. Also based on the within cluster correlation structure Wang [25] proposes an estimator of the regression function, which again has smaller variance than the working independence one. Both Wang's approach and Welsh, Lin and Carroll's rely on the within cluster correlation structure being *known* although in reality it is often unknown. Our estimator $\hat{\Sigma}$ can be used to substitute the (unknown) within cluster correlation structure in their estimation procedure. This would not change the efficiency of the estimator of the regression function because $\hat{\Sigma}$ enjoys convergence rate $n^{-1/2}$. Further, the established asymptotic normality is also useful for statistical inference.

THEOREM 2. *Under the technical conditions in the Appendix, we have*

$$n^{1/2}\bigg\{\hat{\sigma}^2 - \sigma^2 - \frac{1}{4}h^4\bigg(\frac{\mu_1\mu_3 - \mu_2^2}{\mu_0\mu_2 - \mu_1^2}\bigg)^2 b\bigg\}$$
$$\xrightarrow{D} N(0,\ 2\sigma^4(1+\gamma)c_1 + \text{var}(\varepsilon_{11}^2)\gamma c_1).$$

Similarly to Theorem 1, if we choose the bandwidth $h$ to be $O(n^{-1/8})$, the estimator $\hat{\sigma}^2$ will have the convergence rate $n^{-1/2}$.

We give the proofs of these two theorems in the Appendix.

**4. Simulation study.** In this section, we conduct a simulation study on the efficacy of the proposed estimation method, and compare our results with restricted maximum likelihood estimation (REML), which is commonly used in the literature for the estimation of the within cluster correlation structure. We will also demonstrate that the proposed estimator of the covariance matrix can be used to improve the estimators of the functional coefficients, and the proposed estimator of the covariance matrix is robust against likelihood specification and mild model misspecification.

TABLE 1
*The MSEs of the estimators*

|  | 6 | 7 | 8 | 9 | 10 | Our method |
|---|---|---|---|---|---|---|
| $\sigma_{11}$ | 0.1008 | 0.1013 | 0.0997 | 0.0987 | 0.0972 | 0.0967 |
| $\sigma_{12}$ | 0.0785 | 0.0788 | 0.0791 | 0.0789 | 0.0787 | 0.0760 |
| $\sigma_{22}$ | 0.0857 | 0.0874 | 0.0889 | 0.0890 | 0.0899 | 0.0887 |
| $\sigma^2$ | 0.0787 | 0.0461 | 0.0242 | 0.0154 | 0.0095 | 0.0081 |

The top row is the number of knots. The last column is the MSE of the estimators obtained by our method. The rest are the MSE of the estimators obtained by REML at different numbers of knots.

EXAMPLE. In (1.3), with $p=2$, $X_{ij}$ are i.i.d. from $N(0,I_2)$, $U_{ij}$ are i.i.d. from $U(0,1)$ and $\varepsilon_{ij}$ are i.i.d. from $N(0,\sigma^2)$. With $q=2$, $Z_{ij}$ are i.i.d. from $N(0,I_2)$ and $\mathbf{e}_i$ are i.i.d. from $N(0,\Sigma)$. Next, $n_i$ are i.i.d. and set to be the integer part of $|\theta|+6$, $\theta \sim N(0,4)$. We set $m$ equal to 100 and $\sigma^2$ to 1. $\Sigma = (\sigma_{ij})_{2\times 2}$, and $\sigma_{11}$ is set to be 2, $\sigma_{12}$ to 1.5 and $\sigma_{22}$ to 2. We also set $a_1(U) = \sin(2\pi U)$ and $a_2(U) = \cos(2\pi U)$.

The kernel function involved in the local linear modeling is taken to be the Epanechnikov kernel $K(t) = 0.75(1-t^2)_+$. The bandwidth is chosen to be 0.15. We repeat the simulation 100 times; the mean squared error (MSE) is used to assess the accuracy of the estimators. The MSEs of the estimators of $\sigma_{ij}$ and $\sigma^2$ are presented in the last column in Table 1, which suggests that the proposed estimators perform well.

Next, we compare the proposed estimation with REML. For the nonparametric setting, we use the B-spline decomposition to approximate the functional coefficient. The knots in the B-spline decomposition are equally spaced. We choose the range of number of knots where the REML performs best when we use REML to estimate the $\sigma_{ij}$ and $\sigma^2$. The Downhill Simplex approach (Jacoby, Kowalik and Pizzo [16]) is employed for the optimization involved in REML. The MSEs of the obtained estimators are presented in Table 1.

From Table 1, we can see that the newly proposed method and REML have comparable performance. Given the fact that REML is a likelihood-based method fully utilizing the information provided by data, it should be, in theory, the most efficient as long as the likelihood function is correctly specified. However, in practice, REML may not be practicable because the optimization involved in REML can be problematic when the dimension of $\mathbf{e}_i$ is larger than 3. Specifically, the Downhill Simplex method, which served us well previously, has a tendency of failing to converge in such cases. Although we cannot rule out the possibility of better optimization algorithms, the newly proposed estimation does have the considerable practical advantage of yielding a closed-form solution and being optimization-free.

As one referee has rightly pointed out, another advantage of the proposed estimation over REML is that the former does not rely on the likelihood function, so it is robust against likelihood specification. To examine this point, we set the $\varepsilon_{ij}$ as i.i.d. from a uniform distribution with mean 0 and variance $\sigma^2$. The proposed estimation is employed again to estimate the $\sigma_{ij}$ and $\sigma^2$. The MSEs of the obtained estimators are 0.006 for $\hat{\sigma}^2$, 0.090 for $\hat{\sigma}_{11}$, 0.066 for $\hat{\sigma}_{12}$ and 0.091 for $\hat{\sigma}_{22}$. This suggests that the proposed estimation is indeed robust with respect to likelihood specification.

As mentioned in the Introduction, the covariance matrix of random effects serves two purposes. First, it improves the estimator of the functional coefficient. Lin and Carroll [19] have shown that the estimator based on

Table 2
*The Improvement of the estimators*

| 7 | 8 | 9 | 10 | 11 | 12 | 13 | 14 | 15 |
|---|---|---|---|---|---|---|---|---|
| 1.47 | 2.20 | 2.62 | 2.70 | 2.71 | 2.81 | 2.74 | 2.80 | 2.80 |
| 0.10 | 0.19 | 0.30 | 0.45 | 0.67 | 0.94 | 1.28 | 1.59 | 2.00 |

The top row is the number of knots, the second row is the IMP for $a_1(\cdot)$ and the third row is the IMP for $a_2(\cdot)$.

working independence would be the best when kernel smoothing is used in a nonparametric setting. Welsh, Lin and Carroll [26] have also shown the estimator can be improved by the weighted least squares approach when spline smoothing is used. The following is to explore how much improvement can be achieved when we incorporate the proposed estimator of the covariance matrix of random effects in the latter approach.

For any function $g(\cdot)$, if $\hat{g}(\cdot)$ is an estimator of $g(\cdot)$, the mean integrated squared error (MISE) of $\hat{g}(\cdot)$ is defined as

$$\text{MISE} = \int \{\hat{g}(u) - g(u)\}^2 \, du.$$

Let $\text{MISE}_1$ be the MISE of the estimator of the functional coefficient based on working independence, and $\text{MISE}_2$ the MISE based on the weighted least squares approach, incorporating the proposed estimator of the covariance matrix of random effects. We use $\text{IMP} = (\text{MISE}_1 - \text{MISE}_2)/\text{MISE}_2$ to denote the improvement due to the weighted least squares approach.

We have computed IMP when the number of knots in the B-spline decomposition is greater than 7 and less than 15, the choice being made on empirical grounds. In fact, we found that the MISE of the estimators based on either the weighted least squares approach or working independence is much smaller when the number of knots lies in this range than when it lies outside this range. The obtained IMPs are presented in Table 2, which suggests improvement across all cases and by a substantial margin for the larger number of knots.

Finally, another interesting question is whether the proposed estimation still works when the model is misspecified. We have conducted some investigation, leaving a systematic study to the future. We simulated data from the model

$$y_{ij} = x_{ij1}a_1(U_{ij}) + x_{ij2}a_2(U_{ij}) + g_1(z_{ij1})e_{i1} + g_2(z_{ij2})e_{i2} + \varepsilon_{ij},$$

$j = 1, \ldots, n_i$, $i = 1, \ldots, m$. $X_{ij} = (x_{ij1}, x_{ij2})^{\text{T}}$, $Z_{ij} = (z_{ij1}, z_{ij2})^{\text{T}}$, $\mathbf{e}_i = (e_{i1}, e_{i2})^{\text{T}}$, and $U_{ij}$ and $\varepsilon_{ij}$ are simulated in the same way as before. $n_i$ and $m$ are also set in the same way as before. We set $g_1(z) = z + 0.1\sin(z)$ and $g_2(z) = z + 0.1\sin(z)$. We still set $a_1(U) = \sin(2\pi U)$ and $a_2(U) = \cos(2\pi U)$.

Notice that $t_{ij1} = g_1(z_{ij1})$ cannot be treated as a covariate because $g_1(\cdot)$ is treated as *unknown*. The same remark applies to $g_2(z_{ij2})$. So, the model (1.3) is not the true model and we have a misspecified case here.

The proposed estimation is employed again to estimate the $\sigma_{ij}$ and $\sigma^2$. The MSEs of the obtained estimators are 0.006 for $\hat{\sigma}^2$, 0.104 for $\hat{\sigma}_{11}$, 0.076 for $\hat{\sigma}_{12}$ and 0.098 for $\hat{\sigma}_{22}$. This suggests that the proposed estimation is robust against a mild degree of misspecification.

**5. Real data analysis.** The data come from the Bangladesh Demographic and Health Survey (BDHS) of 1996–1997 (Mitra et al. [23]), which is a cross-sectional, nationally representative survey of ever-married women aged between 10 and 49. The analysis is based on a sample of 8189 women nested within 296 primary sampling units or clusters, with sample sizes ranging from 16 to 58. We allow for the hierarchical structure of the data by fitting a two-level model with women at level 1 nested within clusters at level 2. In the multilevel model, cluster-level random effects allow for correlation between outcomes for women in the same cluster. A further hierarchical level is the administrative division; Bangladesh is divided into six administrative divisions which are Barisal, Chittagong, Dhaka, Kulna, Rajshahi and Sylhet. Effects at this level are represented in the model by fixed effects since there are only six divisions.

The dependent variable, $y_{ij}$, is the duration in months between marriage and the first birth for the $j$th woman in the $i$th cluster. A small number of women (0.6% of the total sample size) reported a premarital birth, and these are excluded from the analysis. When a woman was asked for the date of her first marriage in the BDHS, the intention was to collect the age at which she started to live with her husband. However, it is likely that some older women reported the age at which they were formally married, which in Bangladesh can take place at a very young age and some time before puberty (Mitra et al. [23]). For this reason, we calculate the first birth interval assuming a minimum *effective* age at marriage of 12 years. The youngest age at first birth in the sample was 12 years and this is assumed to be the youngest age at which a woman can give birth. 11.53% of women in the sample had not had a birth by the time of the survey and are therefore right-censored.

We consider several covariates which are commonly found to be associated with fertility in Bangladesh. The selected individual-level categorical covariates ($Z_{ij}$) include the woman's level of education (none coded by 0, primary or secondary plus coded by 1), religion (Muslim coded by 1, Hindu or other coded by 0) and age at first marriage in years. Another individual-level covariate is year of marriage ($U_{ij}$). We also consider two cluster-level variables, administrative division and type of region of residence (urban coded by 1, rural coded by 0). We take Barisal as the reference and the differences among the six administrative divisions are modeled by a set of dummy variables.

We take urban as the reference and the differences between urban and rural clusters are modeled by dummy variables. $\xi_i$ is the vector of these six dummy variables.

Typically, there are two ways to analyze right-censored data. One is the likelihood function approach based on Cox proportional hazard function (Cox [4]); another is the regression approach based on an unbiased transformation (Fan and Gijbels [6]). In this paper, we adopt the latter approach. We recover the censored $y_{ij}$ by the unbiased transformation proposed by Fan and Gijbels [6] first, then let $X_{ij} = (Z_{ij}^{\mathrm{T}}, \xi_i^{\mathrm{T}})^{\mathrm{T}}$ and employ (1.3) to fit the transformed data. The proposed estimation procedure is used to estimate the impacts of the covariates as well as the correlations of these impacts. The results obtained are presented in Figure 1 and Table 3.

Table 3 shows that the correlation between the impact of age (of the woman at first marriage) and the impact of education is negative. This implies that the impact of age on the first birth interval would be weak in areas where the education level is high. The impact of age and the impact of religion are strongly negatively correlated. This suggests that the impact of age on the first birth interval is also very weak in areas which are predominantly Muslim. The correlation between the impact of education and the impact of religion is also negative. This suggests that education would not have a big impact on the first birth interval in areas which are predominantly Muslim.

From Figure 1, we can see the trend of length of the first birth interval is decreasing with time. This is attributed to a successful national family planning program (see, e.g., Cleland et al. [2]), which increases the age at first marriage. A nationally representative survey of women in 1996–1997 (Mitra et al. [23]) found that the median age at marriage was 13.3 years among respondents aged 45–49 at the time of the survey, compared to 15.3 years for respondents aged 20–24.

The impact of the woman's age is negative and decreasing with time. The impact of the woman's education is negative until around 1984. Before 1984, the longer birth intervals among women with no education may be partly explained by the higher frequency with which these women report their age at formal marriage rather than their age at cohabitation. Calculating the duration to the first birth from an origin of age 12 for these women may have artificially inflated the lengths of their birth intervals.

Urban impact is negative before 1959, and is getting smaller with time after 1959. This is because at earlier times, some women in rural areas got married very early. There was a considerable time period between their formal marriage and their age at cohabitation. Such cases are getting fewer with time. The impact of following the Muslim religion is always negative, and was decreasing sharply from 1955 to 1968, after which it stayed steady. This suggests that Muslims tend to have significantly shorter first birth

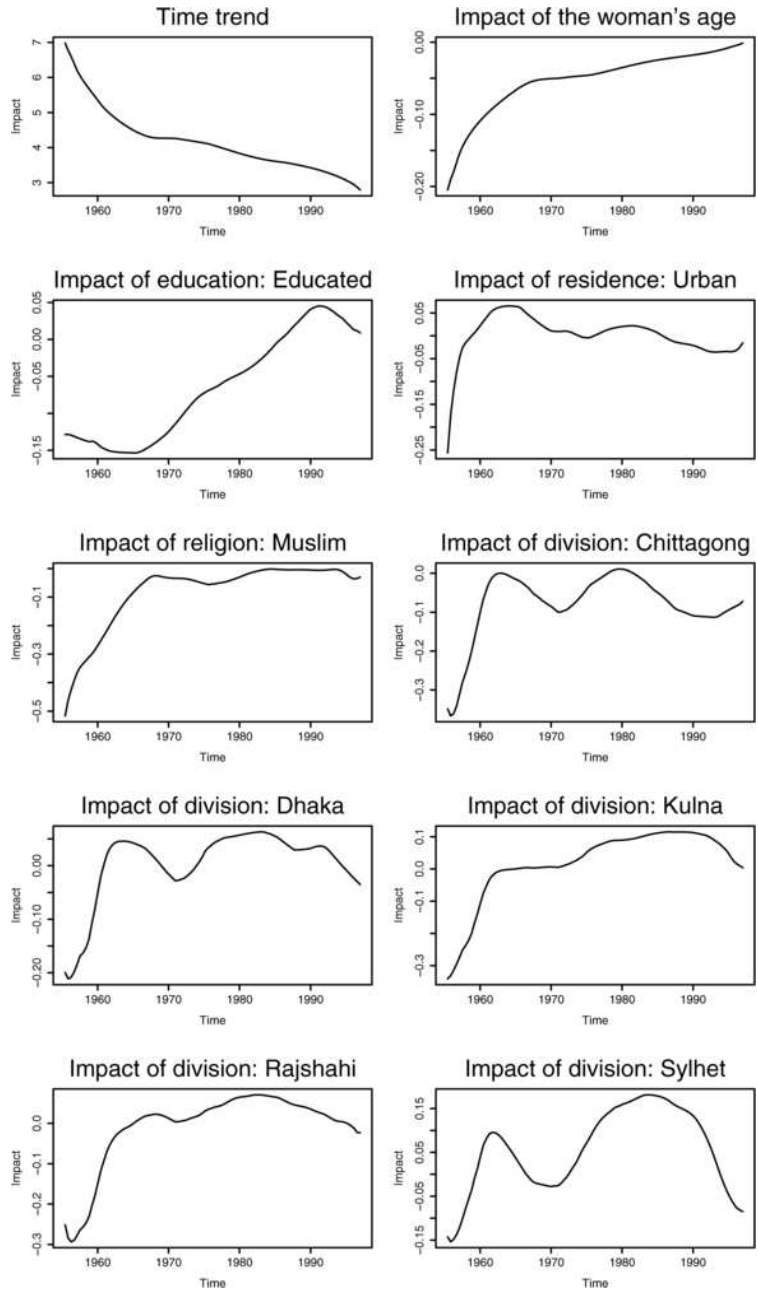

Fig. 1. *The impacts of the covariates which are commonly found to be associated with fertility in Bangladesh on the length of the first birth interval.*

TABLE 3
*The correlation between the impacts of the covariates*

|         | **AGE1MAR** | **EDUC** | **REL** |
|---------|-------------|----------|---------|
| AGE1MAR | 1           | −0.180   | −0.934  |
| EDUC    | −0.180      | 1        | −0.171  |
| REL     | −0.934      | −0.171   | 1       |

AGE1MAR is the impact of the age of the woman at first marriage in years, EDUC is the impact of woman's education and REL is the impact of the woman's religion.

intervals than others before 1968; after 1968, they still tend to have shorter first birth intervals than others but not as significantly.

Looking at the impact of the division, it is noticeable that the intervals are shorter in Chittagong than in the other divisions. This regional effect is as expected and is most likely explained by lower contraceptive use in Chittagong (the most religiously conservative part of Bangladesh) compared to other divisions. Moreover, the impact of the division clearly varies with time.

## APPENDIX

Let $\mathcal{D} = \{(U_{ij}, X_{ij}, Z_{ij}) : j = 1, \ldots, n_i, i = 1, \ldots, m\}$, and we use $(U, X, Z)$ to represent its population. Further, we write $\Omega_1(u) = E(XX^T | U = u)$.

The following technical conditions are imposed to establish the asymptotic results:

(1) $E\varepsilon_{11}^4 < \infty$, $E\|\mathbf{e}_1\|^4 < \infty$, $Ex_i^{2s} < \infty$ and $Ez_j^{2s} < \infty$, where $\|\mathbf{e}_1\| = (\mathbf{e}_1^T \mathbf{e}_1)^{1/2}$, $x_i$ denotes the $i$th element of $X$ and $z_j$ denotes the $j$th element of $Z$ for $s > 2$, $i = 1, \ldots, p$, $j = 1, \ldots, q$.
(2) $a_j''(\cdot)$ is continuous in a neighborhood of $u$ for $j = 1, \ldots, p$, where $a_j''(\cdot)$ is the $j$th element of $\mathbf{a}''(\cdot)$. Further, assume $a_j''(u) \neq 0$.
(3) The marginal density $f(\cdot)$ of $U$ has a continuous derivative in some neighborhood of $u$, and $f(u) \neq 0$.
(4) $r_{ij}(u), \beta_{il}(u)$ and $\gamma_{ij}(u)$ are continuous in a neighborhood of $u$, where

$$r_{ij}(u) = E(x_i x_j | U = u), \qquad \beta_{il}(u) = E(x_i z_l | U = u),$$

$$\gamma_{ij}(u) = E(x_i Z^T \Sigma Z x_j | U = u) \qquad \text{for } i, j = 1, \ldots, p, l = 1, \ldots, q.$$

(5) The function $K(t)$ is a density function with a compact support.
(6) $h \to 0$, $nh^2 \to \infty$ and $nh^8$ is bounded.
(7) There exists a sequence of positive real numbers $M_n$ such that $M_n \to \infty$ and

$$n^{-1} M_n \max_{1 \leq i \leq m, 1 \leq j \leq n_i} \sum_{s=1}^{n_i} (Z_{ij}^T (\mathbf{Z}_i^T \mathbf{Z}_i)^{-1} Z_{is})^2 \xrightarrow{P} 0.$$

For easy reference, we first present some useful lemmas.

LEMMA A.1. *Let*

$$T_{(n)} = \sum_{i=1}^{n} \sum_{j=1}^{n} a_{ij} X_i X_j$$

*be a quadratic form in independent random variables $X_i$ [$E(X_i) = 0$, $E(X_i^2) = 1$], with $\lambda_1, \ldots, \lambda_n$ the eigenvalues of the symmetric matrix $(a_{ij})$, with $a_{ii} = 0$ for all $i$. We denote $\sigma_n^2 = E(T_{(n)}^2)$. Suppose that there is a sequence of real numbers $K(n)$ such that*

(a) $K(n)^2 \sigma_n^{-2} \max_{1 \leq i \leq n}(\sum_{1 \leq j \leq n} a_{ij}^2) \to 0$, $n \to \infty$, *and*

(b) $\max_{1 \leq i \leq n}(E[X_i^2 I_{\{|X_i| > K(n)\}}]) \to 0$, $n \to \infty$, *and that the eigenvalues of the matrix $(a_{ij})$ are negligible:*

(c) $\sigma_n^{-2} \max_{1 \leq i \leq n}(\lambda_i^2) \to 0$, $n \to \infty$; *then $\sigma_n^{-1} T_{(n)}$ has an asymptotic $N(0,1)$ distribution.*

See Commenges and Jacqmin-Gadda ([3], Theorem 1).

LEMMA A.2. *Let $(X_1, Y_1), \ldots, (X_n, Y_n)$ be i.i.d. random vectors, where the $Y_i$'s are scalar random variables. Assume further that $E|Y|^s < \infty$ and $\sup_x \int |y|^s f(x,y)\,dy < \infty$, where $f$ denotes the joint density of $(X,Y)$. Let $K$ be a bounded positive function with bounded support, satisfying a Lipschitz condition. Then*

$$\sup_{x \in D}\left| n^{-1} \sum_{i=1}^{n} \{K_h(X_i - x)Y_i - E[K_h(X_i - x)Y_i]\} \right| = O_P[\{nh/\log(1/h)\}^{-1/2}]$$

*provided that $n^{2\varepsilon - 1} h \to \infty$ for some $\varepsilon < 1 - s^{-1}$.*

This follows immediately from the result obtained by Mack and Silverman [22].

As we need the result of Theorem 2 to prove Theorem 1, we prove Theorem 2 first.

PROOF OF THEOREM 2. On using similar arguments as in Fan and Zhang [11], the asymptotic conditional bias and covariance of $\hat{\mathbf{a}}(U_{ij})$ are equal to

(A.1) $\quad \text{bias}\{\hat{\mathbf{a}}(U_{ij})|\mathcal{D}\} = -\frac{1}{2}h^2 \frac{\mu_1 \mu_3 - \mu_2^2}{\mu_0 \mu_2 - \mu_1^2} \mathbf{a}''(U_{ij})(1 + o_P(1))$

and

(A.2) $\quad \text{cov}\{\hat{\mathbf{a}}(U_{ij})|\mathcal{D}\} = O_P((nh)^{-1}).$

Let $Q_i = I_{n_i} - P_i$, and let $\tilde{Q}_i$ be a diagonal matrix generated from the diagonal elements of $Q_i$. Now

$$
\begin{aligned}
\hat{\sigma}^2 =\ & \frac{1}{n-qm} \sum_{i=1}^m \varepsilon_i^{\mathrm T}(Q_i - \tilde{Q}_i)\varepsilon_i + \frac{1}{n-qm} \sum_{i=1}^m \varepsilon_i^{\mathrm T} \tilde{Q}_i \varepsilon_i \\
& + \frac{1}{n-qm} \sum_{i=1}^m E[(\hat{\mathbf{r}}_i - \mathbf{r}_i)^{\mathrm T}|\mathcal{D}] Q_i E[(\hat{\mathbf{r}}_i - \mathbf{r}_i)|\mathcal{D}] \\
& + \frac{2}{n-qm} \sum_{i=1}^m \{(\hat{\mathbf{r}}_i - \mathbf{r}_i) - E[(\hat{\mathbf{r}}_i - \mathbf{r}_i)|\mathcal{D}]\}^{\mathrm T} Q_i E[(\hat{\mathbf{r}}_i - \mathbf{r}_i)|\mathcal{D}] \\
& + \frac{1}{n-qm} \sum_{i=1}^m \{(\hat{\mathbf{r}}_i - \mathbf{r}_i) - E[(\hat{\mathbf{r}}_i - \mathbf{r}_i)|\mathcal{D}]\}^{\mathrm T} \\
& \qquad\qquad \times Q_i \{(\hat{\mathbf{r}}_i - \mathbf{r}_i) - E[(\hat{\mathbf{r}}_i - \mathbf{r}_i)|\mathcal{D}]\} \\
& + \frac{2}{n-qm} \sum_{i=1}^m \varepsilon_i^{\mathrm T} Q_i E[(\hat{\mathbf{r}}_i - \mathbf{r}_i)|\mathcal{D}] \\
& + \frac{2}{n-qm} \sum_{i=1}^m \varepsilon_i^{\mathrm T} Q_i \{(\hat{\mathbf{r}}_i - \mathbf{r}_i) - E[(\hat{\mathbf{r}}_i - \mathbf{r}_i)|\mathcal{D}]\} \\
\equiv\ & J_{n1} + J_{n2} + J_{n3} + J_{n4} + J_{n5} + J_{n6} + J_{n7}.
\end{aligned}
\tag{A.3}
$$

As $Q_i$ is an idempotent matrix and all the diagonal components of $Q_i - \tilde{Q}_i$ are equal to zero, by straightforward calculation it follows that

$$E(J_{n1}|\mathcal{D}) = \frac{\sigma^2}{n-qm} \sum_{i=1}^m \mathrm{tr}(Q_i - \tilde{Q}_i) = 0,$$

$$E(J_{n2}|\mathcal{D}) = \frac{\sigma^2}{n-qm} \sum_{i=1}^m \mathrm{tr}(\tilde{Q}_i) = \sigma^2,$$

$$\mathrm{cov}(J_{n1}, J_{n2}|\mathcal{D}) = E(J_{n_1} J_{n_2}|\mathcal{D}) = \frac{2\sigma^4}{(n-qm)^2} \sum_{i=1}^m \mathrm{tr}((Q_i - \tilde{Q}_i)\tilde{Q}_i) = 0,$$

$$\mathrm{var}(J_{n1}|\mathcal{D}) = \frac{2\sigma^4}{n-qm} \left\{ \frac{qm - \sum_{i=1}^m \sum_{j=1}^{n_i} (Z_{ij}^{\mathrm T}(\mathbf{Z}_i^{\mathrm T}\mathbf{Z}_i)^{-1} Z_{ij})^2}{n-qm} \right\},$$

$$\mathrm{var}(J_{n2}|\mathcal{D}) = \frac{\mathrm{var}(\varepsilon_{11}^2)}{(n-qm)^2} \sum_{i=1}^m \sum_{j=1}^{n_i} [1 - Z_{ij}^{\mathrm T}(\mathbf{Z}_i^{\mathrm T}\mathbf{Z}_i)^{-1} Z_{ij}]^2.$$

By the law of large numbers, it can be easily verified that

$$n\,\mathrm{var}(J_{n2}|\mathcal{D}) \xrightarrow{P} \mathrm{var}(\varepsilon_{11}^2)\gamma c_1, \qquad n\,\mathrm{var}(J_{n1}|\mathcal{D}) \xrightarrow{P} 2\sigma^4(1+\gamma)c_1.$$

Since the eigenvalues of an idempotent matrix are either 1 or 0, by $E\varepsilon_{11}^4 < \infty$, condition (7) and Lemma A.1, we obtain that

$$n^{1/2} J_{n1} \xrightarrow{D} N(0, 2\sigma^4(1+\gamma)c_1).$$

As $\tilde{Q}_i$ is a diagonal matrix, $J_{n2}$ is a sum of independent variables. By $E\varepsilon_{11}^4 < \infty$ and condition (7), it follows from the Lindeberg–Feller theorem that

$$n^{1/2}(J_{n2} - \sigma^2) \xrightarrow{D} N(0, \mathrm{var}(\varepsilon_{11}^2)\gamma c_1).$$

Since the two terms are uncorrelated, we have that

$$(A.4) \quad n^{1/2}(J_{n1} + J_{n2} - \sigma^2) \xrightarrow{D} N(0, 2\sigma^4(1+\gamma)c_1 + \mathrm{var}(\varepsilon_{11}^2)\gamma c_1).$$

It follows from (A.1) that

$$
\begin{aligned}
J_{n3} &= \frac{1}{n-qm} \sum_{i=1}^{m} \sum_{k=1}^{n_i-q} \sum_{j=1}^{n_i} \sum_{s=1}^{n_i} Q_{ikj} Q_{iks} X_{ij}^{\mathrm{T}} \mathrm{bias}\{\hat{\mathbf{a}}(U_{ij})|\mathcal{D}\} \\
&\qquad\qquad\qquad\qquad\qquad \times X_{is}^{\mathrm{T}} \mathrm{bias}\{\hat{\mathbf{a}}(U_{is})|\mathcal{D}\} \\
&= \frac{1}{4} h^4 \left(\frac{\mu_1\mu_3 - \mu_2^2}{\mu_0\mu_2 - \mu_1^2}\right)^2 b(1 + o_P(1)),
\end{aligned}
$$
(A.5)

where

$$\sum_{l=1}^{n_i} Q_{irl} Q_{isl} = \delta_{rs} = \begin{cases} 1, & r = s, \\ 0, & r \neq s. \end{cases}$$

In the following, we will show that the remaining parts $J_{n4}$ to $J_{n7}$ in (A.3) satisfy $n^{1/2} J_{nl} = o_P(1)$, $l = 4, \ldots, 7$.

By (A.1), (A.2) and the law of large numbers, we have

$$
\begin{aligned}
&E\{|J_{n4}||\mathcal{D}\} \\
&= \frac{h^2|\mu_1\mu_3 - \mu_2^2|}{(n-qm)|\mu_0\mu_2 - \mu_1^2|} \\
&\quad \times \left\{ \sum_{i=1}^{m} \sum_{k=1}^{n_i-q} \sum_{j=1}^{n_i} \sum_{s=1}^{n_i} |Q_{ikj} Q_{iks} X_{is}^{\mathrm{T}} \mathbf{a}''(U_{is})| \right. \\
&\qquad\qquad\qquad\qquad \left. \times E[|(\hat{\mathbf{r}}_{ij} - \mathbf{r}_{ij}) - E[(\hat{\mathbf{r}}_{ij} - \mathbf{r}_{ij})|\mathcal{D}]||\mathcal{D}] \right\} \\
&\quad \times (1 + o_P(1)) \\
&\leq \frac{h^2|\mu_1\mu_3 - \mu_2^2|}{(n-qm)|\mu_0\mu_2 - \mu_1^2|} (1 + o_P(1))
\end{aligned}
$$
(A.6)

$$\times \sum_{i=1}^{m}(n_i-q)\left\{\sum_{s=1}^{n_i}(X_{is}^{\mathrm{T}}\mathbf{a}''(U_{is}))^2\sum_{j=1}^{n_i}X_{ij}^{\mathrm{T}}\mathrm{cov}(\hat{\mathbf{a}}(U_{ij})|\mathcal{D})X_{ij}\right\}^{1/2}$$

$$=O_p((n^{-1}h^3)^{1/2}).$$

By the inequality $ab \leq 2(a^2+b^2)$, $\sum_{l=1}^{n_i}Q_{irl}^2=1$, for $r=1,\ldots,(n_i-q)$, and (A.2), it follows that

$$E\{|J_{n5}||\mathcal{D}\} \leq \frac{4}{n-qm}\sum_{i=1}^{m}\sum_{j=1}^{n_i}(n_i-q)X_{ij}^{\mathrm{T}}\mathrm{cov}(\hat{\mathbf{a}}(U_{ij})|\mathcal{D})X_{ij}$$

(A.7)
$$=O_P((nh)^{-1}).$$

Using (A.1) and the inequality $\boldsymbol{\eta}_i^{\mathrm{T}}Q_i\boldsymbol{\eta}_i \leq \boldsymbol{\eta}_i^{\mathrm{T}}\boldsymbol{\eta}_i$ due to $Q_i$ being an idempotent matrix, we have

$$E(J_{n6}^2|\mathcal{D})$$
$$\leq \frac{h^4\sigma^2}{(n-qm)^2}\left\{\frac{\mu_1\mu_3-\mu_2^2}{\mu_0\mu_2-\mu_1^2}\right\}^2\sum_{i=1}^{m}\sum_{j=1}^{n_i}X_{ij}^{\mathrm{T}}\mathbf{a}''(U_{ij})\mathbf{a}''(U_{ij})^{\mathrm{T}}X_{ij}(1+o_P(1)).$$

Therefore,

(A.8) $$J_{n6}=O_p(n^{-1/2}h^2).$$

As

$$J_{n7}=\frac{2}{n-qm}\sum_{i=1}^{m}\sum_{k=1}^{n_i-q}\sum_{j=1}^{n_i}\sum_{s=1}^{n_i}Q_{ikj}Q_{iks}\{(\hat{\mathbf{r}}_{ij}-\mathbf{r}_{ij})-E[(\hat{\mathbf{r}}_{ij}-\mathbf{r}_{ij})|\mathcal{D}]\}\varepsilon_{is},$$

it follows from straightforward but tedious calculations and Lemma A.2 that

$$\{(\hat{\mathbf{r}}_{ij}-\mathbf{r}_{ij})-E[(\hat{\mathbf{r}}_{ij}-\mathbf{r}_{ij})|\mathcal{D}]\}\varepsilon_{is}$$
$$=\frac{X_{ij}^{\mathrm{T}}\Omega_1^{-1}(U_{ij})}{nhf(U_{ij})(\mu_0\mu_2-\mu_1^2)}$$
$$\times\left\{-h\mu_2\left[\sum_{r=1}^{m}\sum_{l=1}^{n_r}X_{rl}K_h(U_{rl}-U_{ij})(Z_{rl}e_r+\varepsilon_{rl})\varepsilon_{is}\right]\right.$$
$$\left.+\mu_1\left[\sum_{r=1}^{m}\sum_{l=1}^{n_r}X_{rl}(U_{rl}-U_{ij})K_h(U_{rl}-U_{ij})(Z_{rl}e_r+\varepsilon_{rl})\varepsilon_{is}\right]\right\}$$
$$\times(1+o_P(1)).$$

By boundedness of the kernel function, independence of random errors and random effects, we have that $E(J_{n7}^2|\mathcal{D})=O_p((nh)^{-2})$, that is,

(A.9) $$J_{n7}=O_p((nh)^{-1}).$$

Combining (A.3)–(A.9) and condition (6), we obtain that

$$n^{1/2}\left\{\hat{\sigma}^2 - \sigma^2 - \frac{1}{4}h^4\left(\frac{\mu_1\mu_3 - \mu_2^2}{\mu_0\mu_2 - \mu_1^2}\right)^2 b\right\}$$

(A.10)
$$\xrightarrow{D} N(0, 2\sigma^4(1+\gamma)c_1 + \mathrm{var}(\varepsilon_{11}^2)\gamma c_1). \qquad \Box$$

PROOF OF THEOREM 1. Using standard arguments as in the proof of Theorem 2 and (A.1), the conditional bias of $\hat{\Sigma}$ is

$$\mathrm{bias}\{\mathrm{vec}(\hat{\Sigma})|\mathcal{D}\} = \frac{1}{4}h^4\left(\frac{\mu_1\mu_3 - \mu_2^2}{\mu_0\mu_2 - \mu_1^2}\right)^2 \{\mathrm{vec}(B_2) - \mathrm{vec}(B_1)b\}(1 + o_P(1))$$
$$+ O_P((nh)^{-1}),$$

and by straightforward but tedious calculation, the Lindeberg–Feller theorem and condition (7), it follows that

$$n^{1/2}\{\mathrm{vec}(\hat{\Sigma} - \Sigma) - \mathrm{bias}[\mathrm{vec}(\hat{\Sigma})|\mathcal{D}]\} \xrightarrow{D} N(\mathbf{0}, \Delta c_2),$$

where

$$\Delta = E\{\mathbf{e}_1\mathbf{e}_1^T \otimes \mathbf{e}_1\mathbf{e}_1^T\} - \mathrm{vec}(\Sigma)\mathrm{vec}^T(\Sigma) + \sigma^2\{\Sigma \otimes \Gamma + \Gamma \otimes \Sigma + \Delta_1\}$$
$$+ 2\sigma^4\{\Delta_2 - \Delta_3 + c_1/c_2(1+\gamma)\Gamma\} + \mathrm{var}(\varepsilon_{11}^2)\{\Delta_3 + c_1/c_2\gamma\Gamma\}.$$

Therefore, we have

$$n^{1/2}\left\{\mathrm{vech}(\hat{\Sigma} - \Sigma) - \frac{1}{4}h^4\left(\frac{\mu_1\mu_3 - \mu_2^2}{\mu_0\mu_2 - \mu_1^2}\right)^2\{\mathrm{vech}(B_2) - \mathrm{vech}(B_1)b\}\right\}$$
$$\xrightarrow{D} N(\mathbf{0}, (R_q^T R_q)^{-1} R_q^T \Delta R_q (R_q^T R_q)^{-1} c_2). \qquad \Box$$

**Acknowledgments.** We thank the Editor, Associate Editor and referees for their comments, which have helped us to improve the paper significantly.

## REFERENCES

[1] CAI, Z., FAN, J. and LI, R. (2000). Efficient estimation and inferences for varying-coefficient models. *J. Amer. Statist. Assoc.* **95** 888–902. MR1804446
[2] CLELAND, J., PHILLIPS, J. F., AMIN, S. and KAMAL, G. M. (1994). *The Determinants of Reproductive Change in Bangladesh*: *Success in a Challenging Environment*. The World Bank, Washington.
[3] COMMENGES, D. and JACQMIN-GADDA, H. (1997). Generalized score test of homogeneity based on correlated random effects models. *J. Roy. Statist. Soc. Ser. B* **59** 157–171. MR1436561
[4] COX, D. R. (1972). Regression models and life tables (with discussion). *J. Roy. Statist. Soc. Ser. B* **34** 187–220. MR0341758
[5] DIGGLE, P. J., HEAGERTY, P., LIANG, K.-Y. and ZEGER, S. L. (2002). *Analysis of Longitudinal Data*, 2nd ed. Oxford Univ. Press. MR2049007

[6] Fan, J. and Gijbels, I. (1994). Censored regression: Local linear approximations and their applications. *J. Amer. Statist. Assoc.* **89** 560–570. MR1294083

[7] Fan, J., Huang, T. and Li, R. (2007). Analysis of longitudinal data with semiparametric estimation of covariance function. *J. Amer. Statist. Assoc.* **102** 632–641.

[8] Fan, J. and Li, R. (2004). New estimation and model selection procedures for semiparametric modeling in longitudinal data analysis. *J. Amer. Statist. Assoc.* **99** 710–723. MR2090905

[9] Fan, J., Yao, Q. and Cai, Z. (2003). Adaptive varying-coefficient linear models. *J. R. Stat. Soc. Ser. B Stat. Methodol.* **65** 57–80. MR1959093

[10] Fan, J. and Zhang, J.-T. (2000). Two-step estimation of functional linear models with applications to longitudinal data. *J. R. Stat. Soc. Ser. B Stat. Methodol.* **62** 303–322. MR1749541

[11] Fan, J. and Zhang, W. (1999). Statistical estimation in varying coefficient models. *Ann. Statist.* **27** 1491–1518. MR1742497

[12] Hastie, T. J. and Tibshirani, R. J. (1993). Varying-coefficient models (with discussion). *J. Roy. Statist. Soc. Ser. B* **55** 757–796. MR1229881

[13] He, X., Fung, W. K. and Zhu, Z. Y. (2005). Robust estimation in generalized partial linear models for clustered data. *J. Amer. Statist. Assoc.* **100** 1176–1184. MR2236433

[14] Hedeker, D. and Gibbons, R. (1994). A random-effects ordinal regression model for multilevel analysis. *Biometrics* **50** 933–944.

[15] Hoover, D. R., Rice, J. A., Wu, C. O. and Yang, L.-P. (1998). Nonparametric smoothing estimates of time-varying coefficient models with longitudinal data. *Biometrika* **85** 809–822. MR1666699

[16] Jacoby, S. L. S., Kowalik, J. S. and Pizzo, J. T. (1972). *Iterative Methods for Nonlinear Optimization Problems*. Prentice Hall, Englewood Cliffs, NJ. MR0339435

[17] Jennrich, R. I. and Schluchter, M. D. (1986). Unbalanced repeated-measures models with structured covariance matrices. *Biometrics* **42** 805–820. MR0872961

[18] Laird, N. M. and Ware, J. H. (1982). Random-effects models for longitudinal data. *Biometrics* **38** 963–974.

[19] Lin, X. and Carroll, R. J. (2000). Nonparametric function estimation for clustered data when the predictor is measured without/with error. *J. Amer. Statist. Assoc.* **95** 520–534. MR1803170

[20] Lindstrom, M. J. and Bates, D. M. (1990). Nonlinear mixed effects models for repeated measures data. *Biometrics* **46** 673–687. MR1085815

[21] Longford, N. T. (1987). A fast scoring algorithm for maximum likelihood estimation in unbalanced mixed models with nested random effects. *Biometrika* **74** 817–827. MR0919850

[22] Mack, Y. P. and Silverman, B. W. (1982). Weak and strong uniform consistency of kernel regression estimates. *Z. Wahrsch. Verw. Gebiete* **61** 405–415. MR0679685

[23] Mitra, S. N., Al-Sabir, A., Cross, A. R. and Jamil, K. (1997). Bangladesh Demographic and Health Survey 1996–1997. National Institute of Population Research and Training (NIPORT), Mitra and Associates, Dhaka, Bangladesh and DHS Macro International, Inc., Calverton, MD.

[24] Qu, A. and Li, R. (2006). Quadratic inference functions for varying-coefficient models with longitudinal data. *Biometrics* **62** 379–391. MR2227487

[25] Wang, N. (2003). Marginal nonparametric kernel regression accounting for within-subject correlation. *Biometrika* **90** 43–52. MR1966549

[26] WELSH, A. H., LIN, X. and CARROLL, R. J. (2002). Marginal longitudinal nonparametric regression: Locality and efficiency of spline and kernel methods. *J. Amer. Statist. Assoc.* **97** 482–493. MR1941465
[27] WU, H. and LIANG, H. (2004). Backfitting random varying-coefficient models with time-dependent smoothing covariates. *Scand. J. Statist.* **31** 3–19. MR2042595
[28] WU, H. and ZHANG, J.-T. (2002). Local polynomial mixed-effects models for longitudinal data. *J. Amer. Statist. Assoc.* **97** 883–897. MR1941417
[29] XIA, Y. and LI, W. K. (1999). On the estimation and testing of functional-coefficient linear models. *Statist. Sinica* **9** 735–757. MR1711643
[30] ZEGER, S. L. and DIGGLE, P. J. (1994). Semiparametric models for longitudinal data with application to CD4 cell numbers in HIV seroconverters. *Biometrics* **50** 689–699.
[31] ZEGER, S. L., LIANG, K.-Y. and ALBERT, P. S. (1988). Models for longitudinal data: A generalized estimating equation approach. *Biometrics* **44** 1049–1060. MR0980999
[32] ZHANG, W., LEE, S.-Y. and SONG, X. Y. (2002). Local polynomial fitting in semi-varying coefficient model. *J. Multivariate Anal.* **82** 166–188. MR1918619

Y. SUN
SCHOOL OF ECONOMICS
SHANGHAI UNIVERSITY OF FINANCE AND ECONOMICS
P.R. CHINA

W. ZHANG
INSTITUTE OF MATHEMATICS, STATISTICS
AND ACTUARIAL SCIENCE
UNIVERSITY OF KENT
CANTERBURY, KENT CT2 7NF
UNITED KINGDOM
E-MAIL: w.zhang@kent.ac.uk

H. TONG
DEPARTMENT OF STATISTICS
LONDON SCHOOL OF ECONOMICS
HOUGHTON STREET
LONDON WC2A 2AE
UNITED KINGDOM
E-MAIL: h.tong@lse.ac.uk